\documentclass{elsart}
\usepackage{amsmath}
\usepackage{amssymb}

\newcommand{\card}[1]{\ensuremath{\lvert{#1}\rvert}}
\newcommand{\vect}[1]{\ensuremath{\mathbf{#1}}} % vector
\newcommand{\co}[1]{\ensuremath{\overline{#1}}}
\def\median{\mathop{\rm median}\nolimits}

\begin{document}

\begin{frontmatter}

\title{Characterizations of discrete Sugeno integrals as polynomial functions over distributive lattices}

\author{Miguel Couceiro}
\ead{miguel.couceiro[at]uni.lu}
\address{
Mathematics Research Unit, FSTC, University of Luxembourg \\
6, rue Coudenhove-Kalergi, L-1359 Luxembourg, Luxembourg}

\author{Jean-Luc Marichal\thanksref{JLM}}
\ead{jean-luc.marichal[at]uni.lu}
\thanks[JLM]{Corresponding author}
\address{
Mathematics Research Unit, FSTC, University of Luxembourg \\
6, rue Coudenhove-Kalergi, L-1359 Luxembourg, Luxembourg}

\date{Accepted October 7, 2009}

\begin{abstract}
We give several characterizations of discrete Sugeno integrals over bounded distributive lattices, as particular cases of lattice polynomial
functions, that is, functions which can be represented in the language of bounded lattices using variables and constants. We also consider the
subclass of term functions as well as the classes of symmetric polynomial functions and weighted infimum and supremum functions, and present
their characterizations, accordingly. Moreover, we discuss normal form representations of these functions.
\end{abstract}

\end{frontmatter}

% Changes:
% --------
% definition     --> defn
% proposition    --> prop
% theorem        --> thm
% corollary      --> cor
% remark         --> rem
% lemma          --> lem
% example        --> exmp
%
% \begin{proof}  --> \begin{pf*}{Proof.}
% \end{proof}    --> \qed \end{pf*}

%---------------------------------------------------------------------------------------------- Section 1
\section{Introduction}

We are interested in the so-called (discrete) Sugeno integral, which was introduced by Sugeno in \cite{Sug74,Sug77} over real intervals and then
widely investigated in aggregation theory, due to the many applications in fuzzy set theory, data fusion, decision making, pattern recognition,
image analysis, etc. (for general background, see \cite{BelPraCal07,GraMurSug00} and for a recent reference, see \cite{GraMarMesPap09}). In particular, it plays a relevant role as a preference functional in the qualitative frameworks of multicriteria decision making and decision making under uncertainty (for a survey, see \cite{DubMarPraRouSab01}).

Another appealing feature of the Sugeno integral is that, unlike other well-established aggregation functions, it can be defined over ordered domains (not necessarily linearly ordered), where the usual arithmetic operations are not necessarily available. By focusing on the ordered structure of real intervals, Marichal~\cite{Mar01} observed that Sugeno integrals can be seen as particular
lattice polynomial functions, namely those that are idempotent. This fact enables us to naturally extend the original definition of the
Sugeno integral to idempotent polynomial functions over bounded distributive lattices; see \cite{Marc}. (Recall that a lattice polynomial
function is simply a combination of projections and constant functions using the fundamental lattice operations.)

The notion of lattice polynomial function is a natural and well-established concept in lattice theory (see Birkhoff~\cite{Bir67}, Burris and
Sankappanavar~\cite{BurSan81}, Gr\"atzer~\cite{Grae03}, Rudeanu~\cite{Rud01}) and it appears in complete analogy with classical notions such as
that of real polynomial functions. Indeed, just as polynomial functions of several real variables constitute the most basic functions over the
field of real numbers, the lattice polynomial functions can be seen as the most elementary functions defined on an arbitrary bounded lattice. As
first observed by Goodstein~\cite{Goo67}, by assuming distributivity, lattice polynomial functions become exactly those functions which can be
represented in conjunctive and disjunctive normal forms.

%changed
Clearly, not every function on a bounded lattice $L$ can be represented by a lattice polynomial since such a function is necessarily
nondecreasing and not every function $f\colon L^n \to L$ is nondecreasing. But even in the nondecreasing case, many fail to have such a
representation (take, for instance, the real interval $[0,1]$ and consider $f\colon [0,1] \to [0,1]$ given by $f(x)=0$, if $x\leqslant 0.5$, and
$f(x)=1$, otherwise). This fact raises the question: which nondecreasing functions constitute lattice polynomial functions?

%changed
This problem was first considered by Kindermann~\cite{Kind78} who showed that the polynomial functions on a finite lattice are exactly those
monotone functions ``preserving'' (compatible with) all of its tolerances. For descriptions of the latter, see Czedli and Klukovits
\cite{CzeKlu83} and Chajda \cite{Chajda}.

Motivated by the aggregation theory setting, in \cite{CouMar} the authors took a more direct approach to polynomial functions and provided
characterizations of polynomial functions on arbitrary (possibly infinite) bounded distributive lattices both as solutions of certain functional
equations (such as the median decomposition system, see \cite{Marc}) and in terms of necessary and sufficient conditions which have natural
interpretations in the realm of decision making and aggregation theory.

One of the main purposes of the current paper is to reveal overlaps between the theories of lattice functions and aggregation functions, and present useful applications of this framework. To this extent, in this paper we provide an explanatory view of the characterizations given in \cite{CouMar}, for instance, considering questions of independence (irredundancy) of the characterizing properties, and particularize them to the
special case of Sugeno integrals. Moreover, using this connection, we study certain subclasses of polynomial functions which translate into noteworthy subclasses of Sugeno integrals.

%changed
The structure of the article is as follows. In \S{2}, we recall the basic notions and present the preliminary results needed throughout the
paper. Lattice polynomial functions are then considered in \S{3}. We start by discussing representations of polynomial functions in normal form
(such as the classical disjunctive and conjunctive normal forms) as well as the question of (non)uniqueness of these representations. This is
done in \S\ref{sec:WLP-Rep}. In \S\ref{sec:WLP-Prop}, we present the various properties of polynomial functions which were used in \cite{CouMar}
to completely describe the lattice polynomial functions. \S{4} is then devoted to the characterization of the subclass of Sugeno
integrals as well as some other relevant subclasses, namely, of symmetric functions and of weighted infimum and supremum functions.

%---------------------------------------------------------------------------------------------- Section 2
\section{Basic notions and terminology}
%changed this section and following secs according to Gabors comments

In this section we recall some notions and terminology needed throughout this paper. For further background in lattice theory we refer the
reader to, e.g., Davey and Priestley~\cite{DavPri02}, Gr\"atzer~\cite{Grae03}, and Rudeanu~\cite{Rud01}.

A \emph{lattice} is an algebraic structure $\langle L,\wedge, \vee \rangle$ where $L$ is a nonempty set, called \emph{universe}, and where the
two binary operations $\wedge$ and $\vee$ satisfy the commutative, associative, absorption, and idempotent laws. With no danger of
ambiguity, we will denote a lattice by its universe. A lattice $L$ is said to be \emph{distributive} if, for every $a,b,c\in L$,
$$
a\vee (b\wedge c)= (a\vee b)\wedge (a\vee c)\quad\mbox{or, equivalently,}\quad a\wedge (b\vee c)= (a\wedge b)\vee (a\wedge c).
$$

Throughout this paper, we let $L$ denote an arbitrary bounded distributive lattice with least and greatest elements $0$ and $1$, respectively.
For $a,b\in L$, $a\leqslant b$ simply means that $a\wedge b=a$ or, equivalently, $a\vee b=b$.
 A \emph{chain} is simply a lattice such that for every $a,b\in L$ we have $a\leqslant b$ or
$b\leqslant a$. A subset $S$ of a lattice $L$ is said to be \emph{convex} if for every $a,b\in S$ and every $c\in L$ such that $a\leqslant
c\leqslant b$, we have $c\in S$. For any subset $S\subseteq L$, we denote by $\co{S}$ the convex hull of $S$, that is, the smallest convex
subset of $L$ containing $S$. For every $a,b\in S$ such that $a\leqslant b$, we denote by $[a,b]$ the \emph{interval} $[a,b]=\{c\in L: a\leqslant c\leqslant b\}$. For any integer $n\geqslant 1$, we set $[n]=\{1,\ldots,n\}$.

For an arbitrary nonempty set $A$ and a lattice $L$, the set $L^A$ of all functions from $A$ to $L$ constitutes a lattice under the operations
$$
(f\wedge g)(x)=f(x)\wedge g(x) \quad \textrm{ and } \quad (f\vee g)(x)=f(x)\vee g(x),
$$
for every $f,g\in L^A$. In particular, any lattice $L$ induces a lattice structure on the Cartesian product $L^n$, $n\geqslant 1$, by defining
$\wedge$ and $\vee$ componentwise, i.e.,
\begin{eqnarray*}
(a_1,\ldots ,a_n)\wedge (b_1,\ldots ,b_n) &=& (a_1\wedge b_1, \ldots , a_n\wedge b_n), \\
(a_1,\ldots ,a_n)\vee (b_1,\ldots ,b_n) &=& (a_1\vee b_1, \ldots , a_n\vee b_n).
\end{eqnarray*}
We denote the elements of $L$ by lower case letters $a,b,c,\ldots$, and the elements of $L^n$, $n>1$, by bold face letters
$\vect{a},\vect{b},\vect{c},\ldots$. We also use $\vect{0}$ and $\vect{1}$ to denote the least element and greatest element, respectively, of
$L^n$. For $c\in L$ and $\vect{x}=(x_1,\ldots ,x_n)\in L^n$, set
$$
\vect{x}\wedge c=(x_1\wedge c, \ldots , x_n\wedge c) \quad \textrm{and} \quad \vect{x}\vee c=(x_1\vee c, \ldots , x_n\vee c).
$$

The \emph{range} of a function $f\colon L^{n}\rightarrow L$ is defined by $\mathcal{R}_f=\{f(\vect{x}) : \vect{x}\in L^n\}$. The \emph{diagonal section} of $f$ is the function $\delta_f\colon L\to L$ defined by $\delta_f(x)=f(x,\ldots,x)$. A function $f\colon
L^{n}\rightarrow L$ is said to be \emph{nondecreasing} (\emph{in each variable}) if, for every $\vect{a}, \vect{b}\in L^n$ such that
$\vect{a}\leqslant\vect{b}$, we have $f(\vect{a})\leqslant f(\vect{b})$. Note that if $f$ is nondecreasing, then
$\co{\mathcal{R}}_f=[f(\vect{0}),f(\vect{1})]$. A function $f\colon L^n\to L$ is said to be a \emph{$\wedge$-homomorphism} (resp.\ a \emph{$\vee$-homomorphism}) if, for every $\vect{x},\vect{y}\in L^n$, we have $f(\vect{x}\wedge\vect{y})=f(\vect{x})\wedge f(\vect{y})$ (resp.\ $f(\vect{x}\vee\vect{y})=f(\vect{x})\vee f(\vect{y})$).

We finish this section with the notion of lattice functions that we will be interested in hereafter. The class of \emph{lattice polynomial
functions} (or simply, \emph{polynomial functions}) from $L^n$ to $L$ is defined recursively as follows:
\begin{enumerate}
\item[(i)] For each $i\in [n]$ and each $c\in L$, the projection $\vect{x}\mapsto x_i$ and the constant function $\vect{x}\mapsto c$ are
polynomial functions from $L^n$ to $L$.

\item[(ii)] If $f$ and $g$ are polynomial functions from $L^n$ to $L$, then $f\vee g$ and $f\wedge g$ are polynomial functions from $L^n$ to
$L$.

\item[(iii)] Any  polynomial function from $L^n$ to $L$ is obtained by finitely many applications of the rules (i) and (ii).
\end{enumerate}

We refer to those polynomial functions constructed from projections by finitely many applications of (ii) as \emph{lattice term functions} (or
simply, \emph{term functions}). A well-known example of a term function is the ternary median function, which is given by
\begin{eqnarray*}
\median(x_1,x_2,x_3) &=& (x_1 \wedge x_2)\vee (x_2 \wedge x_3)\vee (x_3 \wedge x_1)\\
&=& (x_1 \vee x_2)\wedge (x_2 \vee x_3)\wedge (x_3 \vee x_1).
\end{eqnarray*}

\begin{fact}\label{nondecreasing}
Every polynomial function $f:L^n\to L$ is nondecreasing.
\end{fact}

%changed according to change in terminology
\begin{rem}
Lattice polynomial functions are also referred to as lattice functions (Goodstein~\cite{Goo67}), algebraic functions (Burris and
Sankappanavar~\cite{BurSan81}), or weighted lattice polynomial functions (Marichal~\cite{Marc}), whereas lattice term functions are also
referred to as lattice polynomials (Birkhoff~\cite{Bir67} and Gr\"atzer~\cite{Grae03}).
\end{rem}

%---------------------------------------------------------------------------------------------- Section 3
\section{Lattice polynomial functions}

In this section, we discuss the (non)uniqueness of normal form representations of polynomial functions and present in an explanatory manner
the characterization of these functions given in \cite{CouMar}.

\subsection{Representations of polynomial functions}\label{sec:WLP-Rep}

Goodstein \cite{Goo67} showed that in the case of bounded distributive lattices, polynomial functions are exactly those which allow
representations in disjunctive and conjunctive normal forms (see Proposition~\ref{DNF} below, first appearing in \cite[Lemma 2.2]{Goo67}; see
also Rudeanu~\cite[Chapter~3,\,\S{3}]{Rud01} for a later reference). In this subsection we study such normal form representations of polynomial
functions. We completely describe all possible disjunctive and conjunctive normal form representations of a given polynomial function and
provide necessary and sufficient conditions which guarantee their uniqueness.

\begin{prop}\label{DNF}
Let $f\colon L^{n}\rightarrow L$ be a function. The following conditions are equivalent:
\begin{enumerate}
\item[(i)] $f$ is a polynomial function.

\item[(ii)] There exists $\alpha \colon 2^{[n]}\rightarrow L$ such that
$
f(\vect{x})=\bigvee_{I\subseteq [n]}(\alpha(I)\wedge \bigwedge_{i\in I} x_i).
$
\item[(iii)] There exists $\beta \colon 2^{[n]}\rightarrow L$ such that
$
f(\vect{x})=\bigwedge_{I\subseteq [n]}(\beta(I)\vee \bigvee_{i\in I} x_i).
$
\end{enumerate}
\end{prop}

The expressions given in $(ii)$ and $(iii)$ of Proposition~\ref{DNF} are usually referred to as the \emph{disjunctive normal form} (DNF)
representation and the \emph{conjunctive normal form} (CNF) representation, respectively, of the polynomial function $f$.

\begin{rem}
Proposition~\ref{DNF} can be easily adjusted to term functions by requiring $\alpha$ and $\beta$ to be nonconstant functions from $2^{[n]}$ to
$\{0,1\}$ and satisfying $\alpha(\varnothing)=0$ and $\beta(\varnothing)=1$, respectively.
\end{rem}

The following corollaries belong to the folklore of lattice theory and are immediate consequences of Theorems D and E in \cite{Goo67}.

\begin{cor}\label{Ext1}
Every polynomial function is completely determined by its restriction to $\{0,1\}^n$.
\end{cor}

\begin{cor}\label{Ext2}
A function $g\colon \{0,1\}^n\rightarrow L$ can be extended to a polynomial function $f\colon L^{n}\rightarrow L$ if and only if it is
nondecreasing. In this case, the extension is unique.
\end{cor}

It is easy to see that the DNF and CNF representations of a polynomial function $f\colon L^{n}\rightarrow L$ are not necessarily unique. For
instance, we have $x_1 \vee (x_1 \wedge x_2) = x_1 = x_1\wedge (x_1 \vee x_2)$.

For each $I\subseteq [n]$, let $\vect{e}_I$ be the element of $L^n$ whose $i$th component is $1$, if $i\in I$, and $0$, otherwise. Let $\alpha_f
\colon 2^{[n]}\rightarrow L$ be the function given by $\alpha_f(I)=f(\vect{e}_I)$ and consider the function $\alpha^{*}_f \colon
2^{[n]}\rightarrow L$ defined by
\[
\alpha^{*}_f(I) =
\begin{cases}
\alpha_f(I), & \text{if $\bigvee_{J\varsubsetneq I}\alpha_f(J) < \alpha_f(I)$,} \\
0, & \text{otherwise.}
\end{cases}
\]

Observe that by nondecreasing monotonicity, $\bigvee_{J \varsubsetneq I}\alpha_f(J)\leqslant \alpha_f(I)$ for every $I\subseteq [n]$,  and if
$\bigvee_{J \varsubsetneq I}\alpha_f(J)=\alpha_f(I)$, then
$$
\big(\alpha_f(I)\wedge \bigwedge_{i\in I} x_i\big)\vee \bigvee_{J\varsubsetneq I}\big(\alpha_f(J)\wedge \bigwedge_{i\in J} x_i\big)=\bigvee_{J
\varsubsetneq I}\big(\alpha_f(J)\wedge \bigwedge_{i\in J} x_i\big).
$$
Thus $\alpha_f$ and $\alpha^{*}_f$ give rise to two, possibly distinct, DNF representations of $f$, i.e.,
$$
f(\vect{x})=\bigvee_{I\subseteq [n]}\big(\alpha_f(I)\wedge \bigwedge_{i\in I} x_i\big)=\bigvee_{I\subseteq [n]}\big(\alpha^{*}_f(I)\wedge
\bigwedge_{i\in I} x_i\big).
$$
For each polynomial function $f\colon L^{n}\rightarrow L$, set
$$
\mathrm{DNF}(f)=\Big\{\alpha \in L^{2^{[n]}} : f(\vect{x})=\bigvee_{I\subseteq [n]}\big(\alpha(I)\wedge \bigwedge_{i\in I} x_i\big)\Big\},
$$
and let $\mathrm{A}(f)$ be the set of all those maps $\alpha \in L^{2^{[n]}}$ such that, for every $I\subseteq [n]$,
\begin{itemize}
\item $\alpha (I)\leqslant \alpha_f(I)$,

\item $\bigvee_{J \subseteq I}\alpha (J)= \alpha_f(I)$  whenever $\bigvee_{J \varsubsetneq I}\alpha_f(J)<\alpha_f(I).$
\end{itemize}

\begin{prop}\label{max-min DNF}
For any polynomial function $f\colon L^n\rightarrow L$, we have $\mathrm{DNF}(f)=\mathrm{A}(f)$.
\end{prop}

\begin{pf*}{Proof.}
Let $\alpha\in\mathrm{A}(f)$ and let $I\subseteq [n]$. For any $J\subseteq I$, we have $\alpha^{*}_f(J)\leqslant \bigvee_{K\subseteq I}\alpha(K)$ and hence
$$
f(\vect{e}_I)=\bigvee_{J\subseteq I}\alpha^{*}_f(J)\leqslant \bigvee_{J\subseteq I}\alpha(J)\leqslant \bigvee_{J\subseteq I}\alpha_f(J)=
\alpha_f(I)=f(\vect{e}_I).
$$
By Corollary \ref{Ext1}, $\alpha\in \mathrm{DNF}(f)$ and hence $\mathrm{A}(f)\subseteq \mathrm{DNF}(f)$.

Now let $\alpha \in \mathrm{DNF}(f)$ and let $I\subseteq [n]$. By definition, $\alpha (I)\leqslant \bigvee_{J\subseteq I}\alpha(J) = f(\vect{e}_I)= \alpha_f(I)$. It follows that $\alpha\in \mathrm{A}(f)$ and hence $\mathrm{DNF}(f)\subseteq \mathrm{A}(f)$.\qed
\end{pf*}

Dually, Let $\beta_f \colon 2^{[n]}\rightarrow L$ be the function given by $\beta_f(I)=f(\vect{e}_{[n]\setminus I})$ and consider the function
$\beta^{*}_f \colon 2^{[n]}\rightarrow L$ defined by
\[
\beta^{*}_f(I) =
\begin{cases}
\beta_f(I), & \text{if $\bigwedge_{J \varsubsetneq
I}\beta_f(J)>\beta_f(I)$,} \\
1, & \text{otherwise.}
\end{cases}
\]

As before, both $\beta_f$ and $\beta^{*}_f$ give rise to two, possibly distinct, CNF representations of $f$. For each polynomial function
$f\colon L^{n}\rightarrow L$, set
$$
\mathrm{CNF}(f)=\Big\{\beta \in L^{2^{[n]}} : f(\vect{x})=\bigwedge_{I\subseteq [n]}\big(\beta(I)\vee \bigvee_{i\in I} x_i\big)\Big\}
$$
and let $\mathrm{B}(f)$ be the set of all those maps $\beta \in L^{2^{[n]}}$ such that, for every $I\subseteq [n]$,
\begin{itemize}
\item  $\beta (I)\geqslant \beta_f(I)$,

\item $\bigwedge_{J \subseteq I}\beta (J)= \beta_f(I)$  whenever $\bigwedge_{J \varsubsetneq I}\beta_f(J)>\beta_f(I).$
\end{itemize}

In complete analogy, we have the following result, dual to Proposition \ref{max-min DNF}.
\begin{prop}\label{max-min CNF}
For any polynomial function $f\colon L^n\rightarrow L$, we have $\mathrm{CNF}(f)=\mathrm{B}(f).$
\end{prop}

\begin{rem}
Propositions~\ref{max-min DNF} and \ref{max-min CNF} were established in \cite[\S{3}]{Marc} when $L$ is a chain, and in that case it was shown
that both $\mathrm{A}(f)$ and $\mathrm{B}(f)$ also constitute chains. As it is easy to verify, this property does not hold in the general case
of bounded distributive lattices. For instance, $\mathrm{DNF}(f)$ is closed under $\vee$ but not necessarily under $\wedge$, and dually,
$\mathrm{CNF}(f)$ is closed under $\wedge$ but not necessarily under $\vee$.
\end{rem}

Using Propositions~\ref{max-min DNF} and \ref{max-min CNF}, we obtain the following result which determines in which cases the CNF and DNF
representations are unique.

\begin{cor}\label{cor:Uniqueness}
Let $f\colon L^n \to L$ be a polynomial function. Then $f$ has a unique DNF (resp.\ CNF) representation if and only if for every $I\subseteq
[n]$,
\begin{enumerate}
\item[(i)] $\bigvee_{J\varsubsetneq I}\alpha_f (J)< \alpha_f(I)$ (resp.\ $\bigwedge_{J\varsubsetneq I}\beta_f (J)> \beta_f(I)$), and

\item[(ii)] there is no $b\in L\setminus \{\alpha_f(I)\}$ (resp.\ $ c\in L\setminus \{\beta_f(I)\}$) such that  $\alpha_f(I)=b \vee
\bigvee_{J\varsubsetneq I}\alpha_f (J)$ (resp.\ $ \beta_f(I)=c\wedge \bigwedge_{J\varsubsetneq I}\beta_f (J)$).
\end{enumerate}
\end{cor}

\begin{pf*}{Proof.}
We show that the result holds for the DNF representation; the other claim can be verified dually. Consider a polynomial function $f\colon L^n
\to L$. Using Proposition~\ref{max-min DNF} it is not difficult to verify that conditions $(i)$ and $(ii)$ suffice to guarantee that
$\mathrm{DNF}(f)$ is a singleton. Indeed, let $\alpha $ be a map in $\mathrm{DNF}(f)$ and, for the sake of a contradiction, suppose that there
exists $I\subseteq [n]$ such that $\alpha (I)< \alpha_f(I)$. By condition $(i)$, $\bigvee_{J\varsubsetneq I}\alpha_f (J)< \alpha_f(I)$, and
since there is no $b\in L\setminus \{\alpha_f(I)\}$ such that $\alpha_f(I)=b \vee \bigvee_{J\varsubsetneq I}\alpha_f (J)$, it follows that
$\alpha (I)\vee \bigvee_{J \subsetneq I}\alpha (J)< \alpha_f(I)$. Thus $\alpha \not \in \mathrm{A}(f)$ which contradicts
Proposition~\ref{max-min DNF}.

Now to see that the converse also holds, note that if for some $I\subseteq [n]$, $\bigvee_{J\varsubsetneq I}\alpha_f (J)= \alpha_f(I)$, then
$\alpha^{*}_f \neq \alpha_f$ and thus $f$ does not have a unique DNF representation. So assume that condition $(i)$ holds for every $I\subseteq
[n]$, but suppose that for some $I'\subseteq [n]$ there exists $b\in L\setminus \{\alpha_f(I')\}$ such that  $\alpha_f(I')=b \vee
\bigvee_{J\varsubsetneq I'}\alpha_f (J)$. Consider the map $\alpha \in L^{2^{[n]}}$ given by $\alpha (I)=\alpha_f(I)$, for every $I\neq I'$, and
$\alpha(I')=b$ otherwise. Clearly, $\alpha \neq \alpha_f$ but $\alpha \in \mathrm{A}(f)$, and hence, by Proposition~\ref{max-min DNF}, $\alpha
\in \mathrm{DNF}(f)$. Thus $f$ does not have a unique DNF representation.\qed
\end{pf*}

\begin{rem}
\begin{enumerate}
\item[(i)] Clearly, condition $(ii)$ of Corollary~\ref{cor:Uniqueness} is redundant when $L$ is a chain. However this is not the case when $L$
contains two incomparable elements $a$ and $b$. For instance, by distributivity we have $((a\vee b)\wedge x)\vee a=(b\wedge x)\vee a$, for every
$x\in L$.

\item[(ii)] Note that $\alpha_f$ is the only isotone set function in $\mathrm{DNF}(f)$ and, similarly, $\beta_f$ is the only antitone set
function in $\mathrm{CNF}(f)$.
\end{enumerate}
\end{rem}

\subsection{Properties and characterizations of polynomial functions}\label{sec:WLP-Prop}
\label{properties}

In this subsection, we recall several properties of polynomial functions and discuss the relations between them. Combinations of these
properties yield the various characterizations of polynomial functions presented in \cite{CouMar}. Throughout this subsection, let $S$ be a
nonempty subset of $L$.

We say that a function $f\colon L^{n}\rightarrow L$ is
\begin{itemize}
\item \emph{$\wedge_{S}$-homogeneous} if, for every $\vect{x}\in L^n$ and every $c\in S$, we have $f(\vect{x}\wedge c) = f(\vect{x})\wedge c$.

\item \emph{$\vee_{S}$-homogeneous} if, for every $\vect{x}\in L^n$ and every $c\in S$, we have $f(\vect{x}\vee c) = f(\vect{x})\vee c$.

\item \emph{$S$-idempotent} if, for every $c\in S$, we have $f(c,\ldots,c)=c$.
\end{itemize}
For every integer $m\geqslant 1$, every vector $\vect{x}\in L^m$, and every $f\colon L^{n}\rightarrow L$, we define $\langle\vect{x}\rangle_f\in
L^m$ as the $m$-tuple $\langle\vect{x}\rangle_f=\median(f(\vect{0}),\vect{x},f(\vect{1}))$, where the right-hand side median is taken componentwise.

\begin{prop}\label{Min-Max-Hom-Med}
A function $f\colon L^{n}\rightarrow L$ is $\wedge_{S}$- and $\vee_{S}$-homogeneous if and only if it satisfies $f(\median({r},\vect{x},{s}))= \median(r,f(\vect{x}),s)$ for every $\vect{x}\in L^n$ and every $r,s\in S$. Furthermore, if $f(\vect{0}),f(\vect{1})\in S$ and $f(\vect{0})\leqslant f(\vect{x})\leqslant
f(\vect{1})$, then $f(\vect{x}) = f(\langle\vect{x}\rangle_f)$.
\end{prop}

\begin{lem}\label{lemma:MinMaxIdem}
If $f\colon L^{n}\rightarrow L$ is $\wedge_{S}$- and $\vee_{S}$-homogeneous, then it is $S$-idempotent.
\end{lem}

\begin{pf*}{Proof.}
If $f$ is $\wedge_{S}$- and $\vee_{S}$-homogeneous, then for any $c\in S$,
$$
f(c,\ldots,c)\wedge c = f(c,\ldots,c) = f(c,\ldots,c)\vee c,
$$
and thus $f$ is $S$-idempotent.\qed
\end{pf*}

\begin{prop}[\cite{CouMar}]\label{wlpHomogeneous}
Let $f\colon L^n\to L$ be a nondecreasing function. Then
\begin{enumerate}
\item[(i)] If $f$ is $\wedge_{\co{\mathcal{R}}_f}$- or $\vee_{\co{\mathcal{R}}_f}$-homogeneous, then it is $\co{\mathcal{R}}_f$-idempotent.
In particular, $\mathcal{R}_f=\co{\mathcal{R}}_f=[f(\vect{0}),f(\vect{1})]$ and hence $f$ has a convex range.

\item[(ii)] If $f$ is a polynomial function, then it is $\wedge_{\co{\mathcal{R}}_f}$- and $\vee_{\co{\mathcal{R}}_f}$-homogeneous.
\end{enumerate}
\end{prop}

We say that a function $f\colon L^{n}\rightarrow L$ is
\begin{itemize}
\item \emph{horizontally $\wedge_S$-decomposable} if, for every $\vect{x}\in L^n$ and every $c\in S$, we have $f(\vect{x}) = f(\vect{x}\vee c)\wedge f([\vect{x}]^c)$, where $[\vect{x}]^c$ is the $n$-tuple whose $i$th component is $1$, if $x_i\geqslant c$, and $x_i$, otherwise.

\item \emph{horizontally $\vee_S$-decomposable} if, for every $\vect{x}\in L^n$ and every $c\in S$, we have $f(\vect{x}) = f(\vect{x}\wedge c)\vee f([\vect{x}]_c)$, where $[\vect{x}]_c$ is the $n$-tuple whose $i$th component is $0$, if $x_i\leqslant c$, and $x_i$, otherwise.
\end{itemize}

\begin{rem}\label{rem:8660}
\begin{enumerate}
\item[(i)] In the realm of aggregation, the homogeneity properties can be interpreted as follows: Aggregating components cut at a certain level is the same as cutting at the same level the aggregated value of those components. Moreover, for $S$-idempotent functions $f\colon L^{n}\rightarrow L$, $\wedge_S$- and $\vee_S$-homogeneity can be reformulated as $f(\vect{x}\wedge c) = f(\vect{x})\wedge f(c,\ldots,c)$ and $f(\vect{x}\vee c) = f(\vect{x})\vee f(c,\ldots,c)$, respectively, with $c\in S$, which reveals the ``homomorphic'' nature of $f$.

\item[(ii)] Horizontal $\wedge_S$-decomposability of a function $f\colon L^{n}\rightarrow L$ can be interpreted as follows: For any
$\vect{x}\in L^n$ and any horizontal $\wedge_S$-decomposition of $\vect{x}$ with respect to a level $c\in S$, namely $\vect{x}
= (\vect{x}\vee c)\wedge [\vect{x}]^c$, $f(\vect{x})$ decomposes with respect to $\wedge$. A similar interpretation holds for horizontal $\vee_S$-decomposability.

\item[(iii)] The concepts of $\wedge_S$- and $\vee_S$-homogeneity were used by
Fodor and Roubens \cite{FodRou95} to axiomatize certain classes of aggregation functions in the case when $S=L$ is the real interval $[0,1]$. The concept of horizontal
$\vee_S$-decomposability was introduced, also in the case when $S=L$ is the real interval $[0,1]$, by Benvenuti et al.~\cite{BenMesViv02} as a general
property of the Sugeno integral, and referred to as ``horizontal maxitivity".
\end{enumerate}
\end{rem}

For any $\vect{x}=(x_1,\ldots ,x_n)\in L^{n}$, any $k\in [n]$, and any $c\in L$, set
$$
\vect{x}^{c}_{k} = (x_1,\ldots, x_{k-1},c,x_{k+1},\ldots ,x_n).
$$
We say that a function $f\colon L^{n}\rightarrow L$ is \emph{median decomposable} if, for every $\vect{x}\in L^n$, $f$ satisfies the
\emph{median decomposition system}
\begin{align}\label{MedDecomposition}
 f(\vect{x})=\median\big(f(\vect{x}^{0}_{k}), x_k, f(\vect{x}^{1}_{k})\big)\qquad (k=1,\ldots, n).
\end{align}

\begin{thm}[{\cite[Theorem 17]{Marc}}]\label{MedianSystem}
The solutions of the median decomposition system (\ref{MedDecomposition}) are exactly the polynomial functions from $L^n$ to $L$.
\end{thm}

\begin{rem}
A comparative study of normal form representations of Boolean functions was presented in Couceiro et al.~\cite{CouFolLeh06} where it was shown
that the so-called \emph{median normal form} representation, in which Boolean functions are expressed as repeated applications of the median
function to variables, negated variables, and constants, provides a more efficient representation than the classical conjunctive normal form,
disjunctive normal form and polynomial representations (the latter are also called Zhegalkin polynomial representations due to \cite{Zhe27} or
Reed-Muller polynomial representations due to \cite{Mul54,Ree54}). Even though algorithms for converting the classical CNF, DNF, and polynomial
representations into this median normal form were provided, no hint was given on how to produce median representations, e.g., from truth tables.
In the case of nondecreasing functions, Theorem~\ref{MedianSystem} naturally leads to a recursive procedure for obtaining median representations
of functions independent from the way functions are given. Indeed, by setting an ordering of variables, say, the canonical ordering of
variables, we can repeatedly apply Theorem~\ref{MedianSystem} to the variables of any given function in order to derive a nested formula made of
medians applied to variables and constants. By making use of tools in \cite{CouFolLeh06}, namely the decomposition of any Boolean function as a
nondecreasing function composed with variables and negated variables, this procedure can be easily extended to any Boolean function.
Unfortunately, this approach seems to produce median expressions which are not optimal in the sense of \cite{CouFolLeh06}. To this extent one
needs to find rules to simplify the median expressions produced by the above algorithm. This constitutes an interesting problem for future
research.
\end{rem}

We say that a function $f\colon L^{n}\rightarrow L$ is \emph{strongly idempotent} if, for every $\vect{x}\in L^n$ and every $k\in [n]$, we have
$$
f(x_1,\ldots,x_{k-1},f(\vect{x}),x_{k+1},\ldots,x_n)=f(\vect{x}).
$$
Moreover, we say that a function $f\colon L^{n}\rightarrow L$, $n>1$, has a \emph{componentwise convex range} if, for every $\vect{a}\in L^n$
and every $k\in [n]$, the unary function $f_{\vect{a}}^k\colon L\to L$, given by $f_{\vect{a}}^k(x)=f(\vect{a}_k^x)$, has a convex range. As we are going to observe (see Remark~\ref{rem:main} $(v)$), componentwise range convexity extends continuity to the realm of ordered structures.

Let $f\colon L^n\to L$ be a function and consider the following properties:
\begin{enumerate}
\item[$(\mathbf{P}_{\wedge})$] $f$ is a \emph{componentwise $\wedge$-homomorphism}, that is, $f(\vect{x}_k^{a\wedge b})=f(\vect{x}_k^a)\wedge f(\vect{x}_k^b)$ for all $\vect{x}\in L^n$, $a,b\in L$, and $k\in [n]$.

\item[$(\mathbf{P}_{\vee})$] $f$ is a \emph{componentwise $\vee$-homomorphism}, that is, $f(\vect{x}_k^{a\vee b})=f(\vect{x}_k^a)\vee f(\vect{x}_k^b)$ for all $\vect{x}\in L^n$, $a,b\in L$, and $k\in [n]$.
\end{enumerate}

\begin{fact}\label{fact:pv}
If $f\colon L^n\to L$ satisfies $\mathbf{P}_{\wedge}$ or $\mathbf{P}_{\vee}$, then it is nondecreasing.
\end{fact}

The following result completely characterizes the class of polynomial functions in terms of the properties above.

\begin{thm}[\cite{CouMar}]\label{mainChar}
Let $f\colon L^{n}\rightarrow L$ be a function. The following conditions are equivalent:
\begin{enumerate}
\item[(i)] $f$ is a polynomial function.

\item[(ii)] $f$ is median decomposable.

\item[(iii)] $f$ satisfies $\mathbf{P}_{\wedge}$ and $\mathbf{P}_{\vee}$, and is strongly idempotent, has a convex range and a componentwise convex range.

\item[(iv)] $f$ is nondecreasing, $\wedge_{\co{\mathcal{R}}_f}$-homogeneous, and $\vee_{\co{\mathcal{R}}_f}$-homogeneous.

\item[(v)] $f$ satisfies $\mathbf{P}_{\vee}$, and is $\wedge_{\co{\mathcal{R}}_f}$-homogeneous and horizontally $\vee_{\co{\mathcal{R}}_f}$-decomposable.

\item[(vi)] $f$ satisfies $\mathbf{P}_{\wedge}$, and is horizontally $\wedge_{\co{\mathcal{R}}_f}$-decomposable and $\vee_{\co{\mathcal{R}}_f}$-homogeneous.

\item[(vii)] $f$ satisfies $\mathbf{P}_{\wedge}$ and $\mathbf{P}_{\vee}$, and is $\co{\mathcal{R}}_f$-idempotent, horizontally $\wedge_{\co{\mathcal{R}}_f}$-decomposable, and horizontally $\vee_{\co{\mathcal{R}}_f}$-decomposable.
\end{enumerate}
\end{thm}

By Propositions~\ref{Min-Max-Hom-Med} and \ref{wlpHomogeneous}, any lattice polynomial function $f\colon L^n\to L$ satisfies
$f(\vect{x})=f(\langle\vect{x}\rangle_f)$. Using this fact, horizontal $\wedge_{\co{\mathcal{R}}_f}$-decomposability (resp.\ horizontal $\vee_{\co{\mathcal{R}}_f}$-decomposability) can
be replaced with horizontal $\wedge_L$-decomposability (resp.\ horizontal $\vee_L$-decomposability) in the assertions $(v)$--$(vii)$ of Theorem~\ref{mainChar}.

\begin{rem}\label{rem:main}
\begin{enumerate}
\item[(i)] In the case when $L$ is a chain, the conditions $\mathbf{P}_{\wedge}$ and $\mathbf{P}_{\vee}$ can be replaced in Theorem~\ref{mainChar} by nondecreasing monotonicity. In this case, the properties involved in the characterization given in $(iii)$ do not make use of the lattice operations. Further relaxations to the conditions of Theorem~\ref{mainChar} can be found in \cite{CouMar2}.

\item[(ii)] Consider the conditions given in assertion $(iii)$ of Theorem~\ref{mainChar}. (a) The unary function $f(x)=x^2$ on the real interval $L=[0,1]$ satisfies all these conditions but strong idempotency. (b) Let $L=\{0,a,1\}\times\{0,1\}$, and consider $f\colon L\rightarrow L$ given by $f(11)=f(10)=f(a1)=a1$ and $f(00)=f(a0)=f(01)=01$. This function satisfies all the conditions except $\mathbf{P}_{\wedge}$ and $\mathbf{P}_{\vee}$. (c) Let $L=\{0,a,1\}$ and consider the function $f\colon L^2\to L$ given by
    $$
    f(x_1,x_2)=\begin{cases} 0,& \mbox{if $0\in\{x_1,x_2\}$},\\ a,& \mbox{if $x_1=a$ and $x_2\geqslant a$},\\ 1,& \mbox{if $x_1=1$ and $x_2\geqslant a$}.\end{cases}
    $$
    This function satisfies all the conditions except that it does not have a componentwise convex range.

%Now, let $L=\{0,a,b,1\}$, where $a\vee b=1$ and $a\wedge b=0$, and consider $f\colon L^2\rightarrow L$ given by
%$f(1,1)=1$, $f(0,0)=0$ and $f(x_1,x_2)=a$ otherwise. It is easy to verify that $f$ is nondecreasing, strongly idempotent, and has a
%componentwise convex range, but $b\not \in \mathcal{R}_f$ and thus it does not have a convex range. Similarly, the function $f\colon
%L^2\rightarrow L$ given by
%\[
%f(x_1,x_2) =
%\begin{cases}
%x, & \text{if $x_1=x_2=x$,} \\
%0, & \text{if $x_1=0$ or $x_2=0$,} \\
%1, & \text{otherwise,}
%\end{cases}
%\]
%is nondecreasing, strongly idempotent, and has a convex range, but it does not have a componentwise convex range, e.g., for $\vect{a}=(a,a)$,
%both $f_{\vect{a}}^1$ and $f_{\vect{a}}^2$ do not have convex ranges.

\item[(iii)] Any Boolean function $f\colon \{0,1\}^n\to\{0,1\}$ satisfying $f(\vect{0})\leqslant f(\vect{x})\leqslant f(\vect{1})$ is
$\wedge_{\co{\mathcal{R}}_f}$-homogeneous, $\vee_{\co{\mathcal{R}}_f}$-homogeneous, horizontally $\wedge_{\co{\mathcal{R}}_f}$-decomposable, and horizontally $\vee_{\co{\mathcal{R}}_f}$-decomposable. Moreover, as soon as $n\geqslant 3$, there are such $\co{\mathcal{R}}_f$-idempotent Boolean functions which are
not nondecreasing, thus showing that nondecreasing monotonicity and the properties $\mathbf{P}_{\wedge}$ and $\mathbf{P}_{\vee}$ (which reduce to nondecreasing monotonicity in the case of chains) are necessary in assertions $(iv)$--$(vii)$ of
Theorem~\ref{mainChar}.

\item[(iv)] The set $\co{\mathcal{R}}_f$ cannot be replaced by $\mathcal{R}_f$ in the assertions $(iv)$--$(vii)$ of Theorem~\ref{mainChar}. To
see this, let $L=\{0,1/3,2/3,1\}$ with the canonical ordering of the elements and define $f\colon L \rightarrow L$ by
\[
f(x) =
\begin{cases}
0, & \text{if $x\leqslant 1/3$,} \\
1, & \text{otherwise.}
\end{cases}
\]
Clearly, $f$ is $\wedge_{\mathcal{R}_f}$-homogeneous, $\vee_{\mathcal{R}_f}$-homogeneous, horizontally $\wedge_{\mathcal{R}_f}$-decomposable, and horizontally
$\vee_{\mathcal{R}_f}$-decomposable, but it is not median decomposable.

\item[(v)] In the special case of real interval lattices, i.e., where $L=[a,b]$ for reals $a\leqslant b$, the property of having a convex
range, as well as the property of having a componentwise convex range, are consequences of continuity. More precisely, for nondecreasing
functions $f\colon [a,b]^n\to\mathbb{R}$, being continuous is equivalent to being continuous in each variable, and this latter property is equivalent
to having a componentwise convex range. In fact, since polynomial functions are continuous, the conditions of having a convex range and a
componentwise convex range can be replaced by continuity in $(iii)$ of Theorem \ref{mainChar}. Also, by Proposition~\ref{wlpHomogeneous}, in
each of the conditions $(iv)$--$(vii)$ we can add continuity and replace $\co{\mathcal{R}}_f$ by $\mathcal{R}_f$.
\end{enumerate}
\end{rem}

%---------------------------------------------------------------------------------------------- Section 4
\section{Some particular classes of lattice polynomial functions}

We now consider some important subclasses of polynomial functions, namely, those of (discrete) Sugeno integrals, of symmetric polynomial
functions, and of weighted infimum and supremum functions, and provide their characterizations, accordingly.

\subsection{Discrete Sugeno integrals}

A function $f\colon L^{n}\rightarrow L$ is said to be \emph{idempotent} if it is $L$-idempotent.

\begin{fact}\label{sugeno}
A polynomial function is $\{0,1\}$-idempotent if and only if it is idempotent.
\end{fact}

In \cite[\S{4}]{Marc}, $\{0,1\}$-idempotent polynomial functions are referred to as (\emph{discrete}) \emph{Sugeno integrals}. They coincide exactly with
those functions $\mathcal{S}_{\mu}\colon L^{n}\rightarrow L$ for which there is a fuzzy measure $\mu$ such that
$$
\mathcal{S}_{\mu}(\vect{x})=\bigvee_{I\subseteq [n]}\big(\mu(I)\wedge \bigwedge_{i\in I} x_i\big).
$$
Here, by a \emph{fuzzy measure} $\mu$, we simply mean a set function $\mu \colon 2^{[n]}\rightarrow L$ satisfying $\mu(I)\leqslant \mu(I')$
whenever $I\subseteq I'$, and $\mu(\varnothing)=0$ and $\mu([n])=1$.

\begin{prop}[{\cite[Proposition~12]{Marc}}]\label{prop:sug}
For any polynomial function $f\colon L^{n}\rightarrow L$ there is a fuzzy measure $\mu \colon 2^{[n]}\rightarrow L$ such that $f(\vect{x})=\median(f(\vect{0}),\mathcal{S}_{\mu}(\vect{x}),f(\vect{1})) = \langle\mathcal{S}_{\mu}(\vect{x})\rangle_f$.
\end{prop}

From Lemma~\ref{lemma:MinMaxIdem}, Theorem~\ref{mainChar} $(iv)$, and Fact~\ref{sugeno}, we immediately obtain the following characterization of
the Sugeno integrals, which was previously established in the case of real variables in \cite[Theorem~4.2]{Mar00c}.

\begin{cor}\label{mainChar2}
A function $f\colon L^{n}\rightarrow L$ is a Sugeno integral if and only if it is nondecreasing, $\wedge_{L}$-homogeneous, and $\vee_{L}$-homogeneous.
\end{cor}

Even though Corollary~\ref{mainChar2} can be derived from condition $(iv)$ of Theorem~\ref{mainChar} by simply modifying the two homogeneity
properties, to proceed similarly with conditions $(v)$ and $(vi)$, it is necessary to add the conditions of $\{1\}$-idempotency and
$\{0\}$-idempotency, respectively. To see this, let $L$ be a chain with at least three elements and consider the unary functions $f(x)=x\wedge
d$ and $f'(x)=x\vee d$, where $d\in L\setminus\{0,1\}$. Clearly, $f$ is $\wedge_{L}$-homogeneous and horizontally $\vee_L$-decomposable and $f'$ is $\vee_{L}$-homogeneous and horizontally $\wedge_L$-decomposable. However, neither $f$ nor $f'$ is a Sugeno integral. To see that these additions are
sufficient, just note that $\wedge_L$-homogeneity (resp.\ $\vee_L$-homogeneity) implies $\{0\}$-idempotency (resp.\ $\{1\}$-idempotency).

\begin{thm}\label{mainChar23}
Let $f\colon L^{n}\rightarrow L$ be a function. The following conditions are equivalent:
\begin{enumerate}
\item[(i)] $f$ is a Sugeno integral.

\item[(ii)] $f$ is $\{0,1\}$-idempotent and median decomposable.

\item[(iii)] $f$ satisfies $\mathbf{P}_{\wedge}$ and $\mathbf{P}_{\vee}$, and is strongly idempotent, has range $L$ and a componentwise convex range.

\item[(iv)] $f$ is nondecreasing, $\wedge_{L}$-homogeneous, and $\vee_{L}$-homogeneous.

\item[(v)] $f$ satisfies $\mathbf{P}_{\vee}$, and is $\{1\}$-idempotent, $\wedge_{L}$-homogeneous, and horizontally $\vee_L$-decomposable.

\item[(vi)] $f$ satisfies $\mathbf{P}_{\wedge}$, and is $\{0\}$-idempotent, horizontally $\wedge_L$-decomposable, and $\vee_{L}$-homogeneous.

\item[(vii)] $f$ satisfies $\mathbf{P}_{\wedge}$ and $\mathbf{P}_{\vee}$, and is $L$-idempotent, horizontally $\wedge_L$-decomposable, and horizontally $\vee_L$-decomposable.
\end{enumerate}
\end{thm}

\subsection{Lattice term functions}

A function $f\colon L^{n}\rightarrow L$ is said to be \emph{conservative} if, for every $\vect{x}\in L^n$, we have
\begin{equation}\label{eq:discretiz}
f(\vect{x})\in \{x_1,\ldots,x_n\}.
\end{equation}

\begin{rem}
Condition (\ref{eq:discretiz}) was used in the binary case by Pouzet et
al.~\cite{PouRosSto98}.
\end{rem}

Clearly, every conservative function is idempotent. Furthermore, if $L$ is a chain, then every term function $f\colon L^{n}\rightarrow L$ is
conservative. This is not necessarily true in the general case of bounded distributive lattices. To see this, let $L=\{0,a,b,1\}$, where $a\vee
b=1$ and $a\wedge b=0$, and let $f(x_1,x_2)=x_1\vee x_2.$ Obviously, $f(a,b)\not\in\{a,b\}$.

Thus, in order to characterize term functions in the general case of bounded distributive lattices, we need to relax conservativeness. We say
that a function $f\colon L^{n}\rightarrow L$ is \emph{weakly conservative} if (\ref{eq:discretiz}) holds for every $\vect{x}\in \{0,1\}^n$. By
definition, every weakly conservative function is $\{0,1\}$-idempotent.

\begin{cor}\label{mainChar3}
Let $f\colon L^{n}\rightarrow L$ be a Sugeno integral. Then $f$ is a term function if and only if it is weakly conservative.
\end{cor}

\subsection{Symmetric polynomial functions}

An important property of functions is that of symmetry which basically translates into saying that each argument of the function has the same
``weight'' on the values of the function. Formally, a function $f\colon L^{n}\rightarrow L$ is \emph{symmetric} if, for every permutation
$\sigma$ on $[n]$, we have $f(x_1,\ldots , x_n)=f(x_{\sigma(1)},\ldots , x_{\sigma(n)})$.

Noteworthy examples of symmetric functions are the so-called \emph{order statistic} functions; see for instance Ovchinnikov~\cite{Ovc96}. For
any arity $n\geqslant 1$ and any $k\in [n]$, the $k$th \emph{order statistic} function is the term function $\mathrm{os}_k\colon
L^{n}\rightarrow L$ defined by
$$
\mathrm{os}_k(\vect{x})=\bigvee _{\substack{I\subseteq [n]\\ \card{I}=n-k+1}}\bigwedge _{i\in I}x_i =\bigwedge _{\substack{I\subseteq [n]\\
\card{I}=k}}\bigvee _{i\in I}x_i.
$$
As a matter of convenience, we set $\mathrm{os}_0=0$ and $\mathrm{os}_{n+1}=1$. It is easy to verify that a polynomial function $f\colon L^n\to
L$ is symmetric if and only if the set function $\alpha_f \colon 2^{[n]}\rightarrow L$, given by $\alpha_f(I)=f(\vect{e}_I)$, is
\emph{cardinality based}, that is, for every $I, I'\subseteq [n]$, if $\card{I}=\card{I'}$, then $\alpha_f(I)=\alpha_f(I')$. Thus, by letting
$w_{f}\colon \{0,\ldots,n\}\rightarrow L$ be the isotone function satisfying $\alpha_f(I)=w_{f}(\card{I})$ for every $I\subseteq [n]$, we get
\begin{eqnarray*}
f(\vect{x}) &=& \bigvee_{I\subseteq [n]} \big(\alpha_f(I)\wedge \bigwedge_{i\in I} x_i\big)
~=~ \bigvee^n_{k=0}\bigg(w_{f}(k)\wedge \bigvee_{\substack{I\subseteq [n]\\ \card{I}=k}}\bigwedge _{i\in I}x_i\bigg)\\
&=& \bigvee^n_{k=0}\big(w_{f}(k)\wedge \mathrm{os}_{n-k+1}(\vect{x})\big).
\end{eqnarray*}

We can make dual observations on the CNF representation of symmetric polynomial functions. Indeed, a polynomial function $f\colon L^n\to L$ is
symmetric if and only if the set function $\beta_f \colon 2^{[n]}\rightarrow L$, given by $\beta_f(I)=f(\vect{e}_{[n]\setminus I})$, is
cardinality based. Thus, by letting $v_{f}\colon \{0,\ldots,n\}\rightarrow L$ be the antitone function satisfying $\beta_f(I)=v_{f}(\card{I})$
for every $I\subseteq [n]$, we get
\begin{eqnarray*}
f(\vect{x}) &=& \bigwedge_{I\subseteq [n]}\big(\beta_f(I)\vee \bigvee_{i\in I} x_i\big)
~=~ \bigwedge^n_{k=0}\bigg(v_{f}(k)\vee \bigwedge_{\substack{I\subseteq [n]\\ \card{I}=k}}\bigvee _{i\in I}x_i\bigg)\\
&=& \bigwedge^n_{k=0}\big(v_{f}(k)\vee \mathrm{os}_{k}(\vect{x})\big).
\end{eqnarray*}
Moreover, we have that $v_{f}(i)=w_{f}(n-i)$ for all $i\in \{0,\ldots,n\}$.

These observations are reassembled in the following result characterizing those polynomial functions which are symmetric. The equivalence
between conditions $(i)$, $(ii)$, and $(iii)$ was observed in \cite[Proposition~9]{DukMar08}.

\begin{thm}\label{thm:SymmWLP}
Let $f\colon L^{n}\rightarrow L$ be a polynomial function. Then the following conditions are equivalent:
\begin{enumerate}
\item[(i)]  $f$ is symmetric,

\item[(ii)]  there is a cardinality based function $\alpha \colon 2^{[n]}\rightarrow L$ such that
$$
f(\vect{x})=\bigvee_{I\subseteq [n]}\big(\alpha(I)\wedge \bigwedge_{i\in I} x_i\big),
$$
\item[(iii)] there is an isotone function $w\colon \{0,\ldots,n\}\rightarrow L$ such that
$$
f(\vect{x})= \bigvee^{n}_{k=0}\big(w(n-k)\wedge \mathrm{os}_{k+1}(\vect{x})\big),
$$
\item[(iv)] there is a cardinality based function $\beta \colon 2^{[n]}\rightarrow L$ such that
$$
f(\vect{x})=\bigwedge_{I\subseteq [n]}\big(\beta(I)\vee \bigvee_{i\in I} x_i\big),
$$
\item[(v)] there is an antitone function $v\colon \{0,\ldots,n\}\rightarrow L$ such that
$$
f(\vect{x})= \bigwedge^n_{k=0}\big(v(k)\vee \mathrm{os}_{k}(\vect{x})\big).
$$
\end{enumerate}
\end{thm}

For every $n\geqslant 1$, the \emph{$(2n+1)$-ary median} function is defined by (see \cite[Chapter~IV]{BarMon70})
$$
\median _{2n+1}(\vect{x})=\mathrm{os}_{n+1}(\vect{x})=\bigvee_{\substack{I\subseteq [2n+1]\\ \card{I}=n+1}}\bigwedge _{i\in I}x_i
=\bigwedge_{\substack{I\subseteq [2n+1]\\ \card{I}=n+1}}\bigvee _{i\in I}x_i.
$$

\begin{cor}\label{medianSym}
\begin{enumerate}
\item[(i)] A polynomial function $f\colon L^{n}\rightarrow L$ is symmetric if and only if there exists an isotone (or, equivalently, antitone)
function $w\colon \{0,\ldots,n\}\rightarrow L$ such that
\begin{equation}\label{eq:WLPmedian2}
f(x_1, \ldots ,x_n)= \median _{2n+1}\big(x_1, \ldots ,x_n, w(0),\ldots ,w(n)\big).
\end{equation}

\item[(ii)] A Sugeno integral $f\colon L^{n}\rightarrow L$ is symmetric if and only if there exists an isotone (or, equivalently,
antitone) function $w\colon \{1,\ldots,n-1\}\rightarrow L$ such that
$$
f(x_1, \ldots ,x_n)= \median _{2n-1}\big(x_1, \ldots ,x_n, w(1),\ldots ,w(n-1)\big).
$$

\item[(iii)] A term function is symmetric if and only if it is an order statistic.
\end{enumerate}
\end{cor}

\begin{pf*}{Proof.}
To see that $(i)$ holds, just note that the function $f\colon L^{n}\rightarrow L$ given by (\ref{eq:WLPmedian2}) is a symmetric polynomial
function. The converse follows from $(iii)$ of Theorem~\ref{thm:SymmWLP}. The statements $(ii)$ and $(iii)$ are consequences of $(i)$.\qed
\end{pf*}

\subsection{Weighted infimum and supremum functions}

We say that a function $f\colon L^{n}\rightarrow L$ is a \emph{weighted infimum} function if there are $w_0,w_1,\ldots ,w_n\in L$ such that
\begin{equation}\label{eq:Wmin}
f(\vect{x})= w_0\wedge \bigwedge_{i\in [n]} (w_i\vee x_i).
\end{equation}
Similarly, we say that $f\colon L^{n}\rightarrow L$ is a \emph{weighted supremum} function if there are $v_0,v_1,\ldots ,v_n\in L$ such that
\begin{equation}\label{eq:Wmax}
f(\vect{x})= v_0\vee \bigvee_{i\in [n]} (v_i\wedge x_i).
\end{equation}
Note that the weights $w_0$ and $v_0$ simply constitute the maximum and the minimum possible values of the function, respectively.

\begin{prop}
If $f\colon L^{n}\rightarrow L$ is a $\wedge$-homomorphism or a $\vee$-homomorphism, then it is nondecreasing. Furthermore, $\wedge$-homomorphicity (resp.\ $\vee$-homomorphicity) implies
horizontal $\wedge$-decomposability (resp.\ horizontal $\vee$-decomposability).
\end{prop}

\begin{thm}\label{min-max}
Let $f\colon L^{n}\rightarrow L$ be a polynomial function. Then
\begin{enumerate}
\item[(i)] $f$ is a weighted infimum function if and only if it is a $\wedge$-homomorphism.

\item[(ii)] $f$ is a weighted supremum function if and only if it is a $\vee$-homomorphism.
\end{enumerate}
\end{thm}

\begin{pf*}{Proof.}
We prove $(i)$. The proof of $(ii)$ follows similarly. Suppose first that $f$ is a weighted infimum function. Then, by distributivity we have
that for every $\vect{x},\vect{y}\in L^n$,
\begin{eqnarray*}
f(\vect{x}\wedge \vect{y}) &=& w_0\wedge \bigwedge_{i\in [n]} \big(w_i\vee(x_i \wedge y_i)\big) ~=~ w_0\wedge \bigwedge_{i\in [n]} \big((w_i\vee x_i)\wedge (w_i\vee y_i)\big)\\
&=& \Big(w_0\wedge \bigwedge_{i\in [n]} (w_i\vee x_i)\Big) \wedge \Big(w_0\wedge \bigwedge_{i\in [n]} (w_i\vee x_i)\Big)\\
&=& f(\vect{x})\wedge f(\vect{y}).
\end{eqnarray*}
In other words, $f$ is a $\wedge$-homomorphism.

Now we show that if $f$ is a $\wedge$-homomorphism, then it is a weighted infimum function. Observe first that, as a polynomial function, $f$ is
$\wedge_{\co{\mathcal{R}}_f}$- and $\vee_{\co{\mathcal{R}}_f}$-homogeneous, and thus, by Proposition~\ref{Min-Max-Hom-Med}, for every $i\in
[n]$ and every $c\in L$,
$$
f(\vect{1}^0_{i}\vee c) = f\big(\langle\vect{1}^0_{i}\vee c\rangle_f\big) = f(\langle\vect{1}^0_{i}\rangle_f)\vee \langle c\rangle_f =
f(\vect{1}^0_{i})\vee \langle c\rangle_f.
$$
Therefore, by $\wedge$-homomorphicity we have
\begin{eqnarray*}
f(\vect{x})&=& f\Big(\bigwedge_{i\in [n]} (\vect{1}^0_{i}\vee x_i)\Big) ~=~ \bigwedge_{i\in [n]} \big(f(\vect{1}^0_{i})\vee \langle x_i\rangle_f\big)\\
&=& \bigwedge_{i\in [n]} \big(f(\vect{1}^0_{i})\vee (f(\vect{1})\wedge(x_i\vee f(\vect{0})))\big)\\
&=& f(\vect{1})\wedge \bigwedge_{i\in [n]} \big(f(\vect{1}^0_{i})\vee x_i\big).
\end{eqnarray*}
Setting $w_0=f(\vect{1})$ and $w_i=f(\vect{1}^0_{i})$ for $i\in [n]$, we have that $f$ is a weighted infimum function.\qed
\end{pf*}

\begin{rem}
\begin{enumerate}
\item[(i)] Idempotent weighted infimum functions $f\colon L^{n}\rightarrow L$ are those functions (\ref{eq:Wmin}) for which $w_0=1$ and
$\wedge_{i\in [n]} w_i=0$. Dually, idempotent weighted supremum functions $f\colon L^{n}\rightarrow L$ are those functions (\ref{eq:Wmax}) for
which $v_0=0$ and $\vee_{i\in [n]} v_i=1$. These functions were introduced on real intervals by Dubois and Prade~\cite{DubPra86} in fuzzy set
theory, and referred to as ``weighted minimum and maximum functions".

\item[(ii)] The restriction of Theorem~\ref{min-max} to idempotent functions (Sugeno integrals) was already established in the special
case of real intervals in \cite[Theorem~5.2]{Mar00c}.
\end{enumerate}
\end{rem}

%---------------------------------------------------------------------------------------------- Section 5
\section{Conclusion and future work}

In this paper we have discussed discrete Sugeno integrals when regarded as particular lattice polynomial functions, namely idempotent polynomial functions. Having this definition in mind, we have investigated normal form representations of polynomial functions and provided complete descriptions of the set of all disjunctive and conjunctive normal forms of a given function, and obtained necessary and sufficient conditions which guarantee uniqueness of such representations. Then, using the characterizations of polynomial functions given in terms of well-established properties in aggregation theory (here generalized to ordered domains where incomparabilities may occur), we presented various axiomatizations for the class of Sugeno integrals over bounded distributive lattices. Furthermore, we discussed redundancy issues concerning the set of axioms and showed that their independence in most cases. Regarding the axiomatization given in Theorem~\ref{mainChar} $(iii)$, it remains open whether or not range convexity is redundant under the other conditions.

Some noteworthy subclasses of Sugeno integrals were also considered, namely those of
\begin{itemize}
\item term functions, which were shown to constitute weakly conservative Sugeno integrals,

\item symmetric Sugeno integrals, which, as in the case of real closed intervals, were shown to constitute Sugeno integrals determined by cardinality based capacities or in fact median functions up to fixing arguments, and

\item weighted infimum and weighted supremum functions, which naturally extend the so-called weighted minimum and weighted maximum functions and which were shown to constitute $\wedge$- and $\vee$-semilattice homomorphisms.
\end{itemize}

Concerning directions of further research, the authors have proposed within this framework, a very natural generalization to the Sugeno integral, namely that of \emph{quasi-Sugeno integral} \cite{CouMar3} and which falls within the scope of utility based decision making. More precisely, these are functions $f\colon L^n\to L$ which can be factorized as a composition
$$
f(x_1,\ldots,x_n)=p(\varphi(x_1),\ldots,\varphi(x_n)),
$$
where $p\colon L^n\to L$ is a Sugeno integral (or equivalently, a polynomial function) and $\varphi\colon L\to L$ a nondecreasing function (utility function). There are several extensions available within this framework, for instance, one could consider possibly different domain and codomain bounded distributive lattices $X$ and $Y$, and remove the monotonicity condition on the inner function $\varphi\colon X\to Y$.

%\bibliographystyle{abbrv}   % styles: plain, unsrt, alpha, abbrv, ieeetr, acm, siam, apalike, amsplain,...
%\bibliography{References}

\end{document}